\documentclass[a4paper,12pt]{article}
\usepackage{amsfonts}

\usepackage{mathrsfs}
\usepackage{algorithm, algpseudocode}
\usepackage{amsmath,amssymb,amsthm,latexsym,amsfonts}
\usepackage{pstricks}
\usepackage{enumerate}

\usepackage{epsfig}
\usepackage{epstopdf}

\usepackage{graphicx}
\usepackage{color}
\usepackage{ifpdf}

\usepackage{caption}

\begin{document}

\setlength{\parindent}{1em}
\renewcommand{\baselinestretch}{1}

\newtheorem{thm}{Theorem}
\newtheorem{definition} {Definition}
\newtheorem{prob}{Problem}
\newtheorem{lem}{Lemma}
\newtheorem{pro}{Proposition}
\newtheorem{cor}{Corollary}
\newtheorem{fact}{Fact}

\theoremstyle{definition}
\newtheorem{defn}{Definition}
\newtheorem{example}{Example}

\def\into{\rightarrow}
\def\is{\leftarrow}
\def\NlBar{\overline{N^\ell}}
\def\-{\mbox{--}}
\newcommand{\bound}[2]{\sum_{i=1}^{#1}{\min\{2i+1, #2\}}}
\newcommand{\halff}[1]{\Big\lfloor \frac{#1}{2} \Big\rfloor}
\newcommand{\halfc}[1]{\Big\lceil  \frac{#1}{2} \Big\rceil}

\newtheorem{claim}{Claim}
\newtheorem{property}{Property}
\newtheorem{remark}{Remark}
\newtheorem{obser}{Observation}
\newtheorem{conj}{Conjecture}
\newtheorem{case}{Case}
\newtheorem{ex}{Example}

\def\pf{\noindent {\it Proof.} }

\title{\bf Erd\H{o}s-Gallai-type results for total monochromatic connection of graphs\footnote{Supported by NSFC No.11371205 and 11531011, and PCSIRT.} }
\author{\small Hui Jiang, Xueliang Li, Yingying Zhang\\
\small Center for Combinatorics and LPMC\\
\small Nankai University, Tianjin 300071, China\\
\small E-mail: jhuink@163.com; lxl@nankai.edu.cn;\\
\small zyydlwyx@163.com}
\date{}
\maketitle
\begin{abstract}

A graph is said to be {\it total-colored} if all the edges and the vertices of the graph are colored. A total-coloring of a graph is a {\it total monochromatically-connecting coloring} ({\it TMC-coloring}, for short) if any two vertices of the graph are connected by a path whose edges and internal vertices have the same color. For a connected graph $G$, the {\it total monochromatic connection number}, denoted by $tmc(G)$, is defined as the maximum number of colors used in a TMC-coloring of $G$. In this paper, we study two kinds of Erd\H{o}s-Gallai-type problems for $tmc(G)$ and completely solve them.


\end{abstract}

\section{Introduction}
In this paper, all graphs are simple, finite and undirected. We refer to book \cite{B} for undefined notation and terminology in graph theory. Throughout this paper, let $n$ and $m$ denote the order (number of vertices) and size (number of edges) of a graph, respectively. Moreover, a vertex of a connected graph is called a {\it leaf} if its degree is one; otherwise, it is an {\it internal vertex}. Let $l(T)$ and $q(T)$ denote the number of leaves and the number of internal vertices of a tree $T$, respectively, and let $l(G)=\max\{ l(T) | $ $T$ is a spanning tree of $G$ $\}$ and $q(G)=\min\{ q(T) | $ $T$ is a spanning tree of $G$ $\}$ for a connected graph $G$. Note that the sum of $l(G)$ and $q(G)$ is $n$ for any connected graph $G$ of order $n$. A path in an edge-colored graph is a {\it monochromatic path} if all the edges on the path have the same color. An edge-coloring of a connected graph is a {\it monochromatically-connecting coloring} ({\it MC-coloring}, for short) if any two vertices of the graph are connected by a monochromatic path of the graph. For a connected graph $G$, the {\it monochromatic connection number} of $G$, denoted by $mc(G)$, is defined as the maximum number of colors used in an MC-coloring of $G$. An {\it extremal MC-coloring} is an MC-coloring that uses $mc(G)$ colors. Note that $mc(G)=m$ if and only if $G$ is a complete graph. The concept of $mc(G)$ was first introduced by Caro and Yuster \cite{CY} and has been well-studied recently. We refer the reader to \cite{CLW,GLQZ} for more details.

In \cite{JLZ}, we introduced the concept of total monochromatic connection of graphs. A graph is said to be {\it total-colored} if all the edges and the vertices of the graph are colored. A path in a total-colored graph is a {\it total monochromatic path} if all the edges and internal vertices on the path have the same color. A total-coloring of a graph is a {\it total monochromatically-connecting coloring} ({\it TMC-coloring}, for short) if any two vertices of the graph are connected by a total monochromatic path of the graph. For a connected graph $G$, the {\it total monochromatic connection number}, denoted by $tmc(G)$, is defined as the maximum number of colors used in a TMC-coloring of $G$. An {\it extremal TMC-coloring} is a TMC-coloring that uses $tmc(G)$ colors. It is easy to check that $tmc(G)=m+n$ if and only if $G$ is a complete graph. Moreover, in \cite{JLZ0} we determined the threshold function for a random graph to have $tmc(G)\geq f(n)$, where $f(n)$ is a function satisfying $1\leq f(n)<\frac{1}{2}n(n-1)+n$. Actually, these concepts are not only inspired by the concept of monochromatic connection number but also by the concepts of monochromatic vertex connection number and total rainbow connection number of connected graphs. For details about them we refer to \cite{CLW1,JLZ1,LMS1,S}. From the definition of the total monochromatic connection number, the following results are immediate.

\begin{pro}\label{pro1}\cite{JLZ} If $G$ is a connected graph and $H$ is a connected spanning subgraph of $G$, then $tmc(G)\geq m(G)-m(H)+tmc(H)$.
\end{pro}

\begin{thm}\label{thm1}\cite{JLZ} For a connected graph $G$, $tmc(G)\geq m-n+2+l(G)$.
\end{thm}

In particular, $tmc(G)=m-n+2+l(G)$ if $G$ is a tree. In \cite{JLZ} we also showed that there are dense graphs that still meet this lower bound.

\begin{thm}\label{thm2}\cite{JLZ} Let $G$ be a connected graph of order $n>3$. If $G$ satisfies any of the following properties, then $tmc(G)=m-n+2+l(G)$.

$(a)$ The complement $\overline{G}$ of $G$ is $4$-connected.

$(b)$ $G$ is $K_3$-free.

$(c)$ $\Delta(G)<n-\frac{2m-3(n-1)}{n-3}$.

$(d)$ $diam(G)\geq 3$.

$(e)$ $G$ has a cut vertex.
\end{thm}

Moreover, we gave an example in \cite{JLZ} to show that the lower bound $m-n+2+l(G)$ is not always attained.

\begin{lem}\label{lem1}\cite{JLZ} Let $G= K_{n_1,\ldots,n_r}$ be a complete multipartite graph with $n_1 \geq \ldots \geq n_t\geq 2$ and $n_{t+1}=\ldots=n_r=1$. Then $tmc(G)=m+r-t$.
\end{lem}
Let $G$ be a connected graph and $f$ be an extremal TMC-coloring of $G$ that uses a given color $c$. Note that the subgraph $H$ formed by the edges and vertices with color $c$ is a tree where the color of each internal vertex is $c$; see \cite{JLZ}. Now we define the {\it color tree} as the tree formed by the edges and vertices with color $c$, denoted by $T_c$. If $T_c$ has at least two edges, the color $c$ is called {\it nontrivial}; otherwise, $c$ is {\it trivial}. We call an extremal TMC-coloring {\it simple} if for any two nontrivial colors $c$ and $d$, the corresponding trees $T_c$ and $T_d$ intersect in at most one vertex. If $f$ is simple, then the leaves of $T_c$ must have distinct colors different from color $c$. Moreover, a nontrivial color tree of $f$ with $m'$ edges and $q'$ internal vertices is said to {\it waste} $m'-1+q'$ colors since the edges and internal vertices of a nontrivial color tree must have the same color. In fact, we can use at most $m+n$ colors to assign its edges and vertices with different colors. Thus, if $f$ wastes $x$ colors, then $tmc(G)=m+n-x$. For the rest of this paper we will use these facts without further mentioning them. In addition, we list a helpful lemma below.

\begin{lem}\label{lem2}\cite{JLZ} Every connected graph $G$ has a simple extremal TMC-coloring.
\end{lem}

Among many interesting problems in extremal graph theory is the Erd\H{o}s-Gallai-type problem to determine the maximum or minimum value of a graph parameter with some given properties. In \cite{CLW,CLW1}, the authors investigated two kinds of Erd\H{o}s-Gallai-type problems for monochromatic connection number and monochromatic vertex connection number, respectively. Motivated by these, we study two kinds of Erd\H{o}s-Gallai-type problems for $tmc(G)$ in this paper.

\noindent$\displaystyle$\textbf{Problem A.} Given two positive integers $n$ and $k$, compute the minimum integer $f(n,k)$ such that for any  connected graph $G$ of order $n$, if $|E(G)|\geq f(n,k)$ then $tmc(G)\geq k$.

\noindent$\displaystyle$\textbf{Problem B.} Given two positive integers $n$ and $k$, compute the maximum integer $g(n,k)$ such that for any  connected graph $G$ of order $n$, if $|E(G)|\leq g(n,k)$ then $tmc(G)\leq k$.

Note that for a connected graph $G$ we have $3\leq tmc(G)\leq\binom{n}{2}+n$, and that $g(n,k)$ does not exist for $3\leq k\leq n-1$ since for a star $S_n$ on $n$ vertices we have $tmc(S_n)=n$. Thus we just need to determine the exact values of $f(n,k)$ for $3\leq k\leq \binom{n}{2}+n$ and $g(n,k)$ for $n\leq k\leq \binom{n}{2}+n$ in the following.

\begin{thm}\label{thm3} Given two positive integers $n$ and $k$ with $3\leq k\leq\binom{n}{2}+n$,
\begin{eqnarray}f(n,k)=
\begin{cases}
n-1 &if\ k=3, \cr
n+k-t-2 &if\ k=\binom{t}{2}+t+2-s,\ where\ 0\leq s\leq t-1 \ and\ 2\leq t\leq n-2, \cr
k &if\ \binom{n}{2}-n+4\leq k\leq \binom{n}{2}+n-3\lfloor\frac{n}{2}\rfloor\ except\ for\ n\
is\ odd\cr
\ \ &\ \ \ \ and\ k=\binom{n}{2}+n-3\lfloor\frac{n}{2}\rfloor,\cr
\binom{n}{2}-r &if\ \binom{n}{2}+n-3(r+1)<k\leq\binom{n}{2}+n-3r,\ where\ 0\leq r\leq \lfloor\frac{n}{2}\rfloor-1\cr
\ \ &\ \ \ \ or\ n\
is\ odd,\ r=\lfloor\frac{n}{2}\rfloor\ and\ k=\binom{n}{2}+n-3\lfloor\frac{n}{2}\rfloor.
 \end{cases}
\end{eqnarray}
\end{thm}

\begin{thm}\label{thm4} Given two positive integers $n$ and $k$ with $n\leq k\leq \binom{n}{2}+n$,
\begin{eqnarray}g(n,k)=
\begin{cases}
k-n+t &if\ \binom{n-t}{2}+t(n-t-1)+n\leq k\leq \binom{n-t}{2}+t(n-t)+n-2, \cr
k-n+t-1 &if\ k=\binom{n-t}{2}+t(n-t)+n-1, \cr
\binom{n}{2}-1 &if\ k=\binom{n}{2}+n-1, \cr
\binom{n}{2} &if\ k=\binom{n}{2}+n,
\end{cases}
\end{eqnarray}
for $2\leq t\leq n-1$.
\end{thm}
In the next sections we will give the proofs of the two theorems. 

\section{Proof of Theorem \ref{thm3}}

Firstly, we give some useful lemmas.

\begin{lem}\label{lem3}\cite{DJS} Let $G$ be a connected graph with $|E(G)|\geq |V(G)|+\binom{t}{2}$ and $t\leq |V(G)|-3$. Then $G$
has a spanning tree with at least $t+1$ leaves and this is best possible.
\end{lem}

Given three nonnegative integers $n$, $t$ and $s$ such that $2\leq t\leq n-2$ and $0\leq s\leq t-1$. We can find a graph $G_{t,s}$ on $n$ vertices with $m(G_{t,s})=n+\binom{t}{2}-1-s$ and $l(G_{t,s})=t$. Construct $G_{t,s}$ as follows: first let $H$ be the graph obtained from a complete graph $K_{t+1}$ by replacing its one edge $uv$ by a path of $n-t$ edges between the ends of $uv$; second we delete $s$ edges between $u$ and the vertices of $V(K_{t+1})\backslash\{u,v\}$ from $H$. It can be checked that $m(G_{t,s})=n+\binom{t}{2}-1-s$ and $l(G_{t,s})=t$. Next we will show that $tmc(G_{t,s})=m(G_{t,s})-n+2+l(G_{t,s})$.

\begin{lem}\label{lem0} Let $G$ be the graph $G_{t,s}$ described above. Then $tmc(G)=m-n+2+l(G).$
\end{lem}
\pf Let $f$ be a simple extremal TMC-coloring of $G$. Suppose that $f$ consists of $k$ nontrivial color trees, denoted by $T_1,\ldots,T_k$. Observe that every vertex appears in at least one of the nontrivial color trees. Suppose $k\geq 2$. Let $T'$ and $T''$ be any two nontrivial color trees of $f$. Since $f$ is simple, there is at most one common vertex between $T'$ and $T''$.
If $T'$ and $T''$ have no common vertex, then there is a total monochromatic path between each vertex of $V(T')$ and each vertex of $V(T'')$.
If $T'$ and $T''$ have a common vertex $w'$, then there is a total monochromatic path between each vertex of $V(T')$ and each vertex of $V(T'')\backslash \{w'\}$. Moreover, $\delta(G)\leq 2$. Hence, $k=2$ and there exists a common vertex $w$ between $T_1$ and $T_2$, which is a leaf of $T_1$ and $T_2$, respectively. In addition, $w$ is the unique vertex of degree two in $G$. If $t<n-2$, then there exist at least two vertices of degree two in $G$, a contradiction. If $t=n-2$, then there exists an edge between the two neighbors of $w$, a contradiction to the construction of $G$. Hence, $k=1$ and so $tmc(G)=m-n+2+l(G)$.
\qed

Given two positive integers $n$ and $p$ with $\frac{n}{2}<p<n-2$, let $t=2(p+1)-n$ and $G_n^{t}$ be the graph defined as follows: partition the vertex set of the complete graph $K_n$ into $n-p$ vertex-classes $V_1,V_2,...,V_{n-p}$, where $|V_1|=|V_2|=...=|V_{n-p-1}|=2$ and $|V_{n-p}|=t$; for each $j\in\{1,...,n-p\}$, select a vertex $v_j^*$ from $V_j$, and delete all the edges joining $v_j^*$ to the other vertices in $V_j$. Next we will show that $tmc(G_n^t)=m(G_n^t)$.

\begin{lem}\label{lem4}Let $G$ be the graph $G_n^t$ described above. Then $tmc(G)=m$.
\end{lem}
\pf Let $f$ be a simple extremal TMC-coloring of $G$. Suppose that $f$ consists of $k$ nontrivial color trees, denoted by $T_1,\ldots,T_k$, where $t_i=|V(T_i)|$ and $q_i=q(T_i)$ for $1\leq i\leq k$. Observe that every vertex appears in at least one of the nontrivial color trees. Note that $m-n+2+l(G)=m$ and $tmc(G)\geq m$ by Theorem \ref{thm1}. As $T_i$ has $t_i-1$ edges and $q_i$ internal vertices, it wastes $t_i-2+q_i$ colors. To show $tmc(G)\leq m$, we just need to show that $f$ wastes at least $n$ colors, i.e. $\sum_{i=1}^k(t_i-2+q_i)\geq n$. In fact, consider the spanning subgraph $G'$ consisting of the union of the $T_i$'s and let $C_1,\ldots,C_s$ denote its components. We claim that for $1\leq i \leq n-p$, the vertices of $V_i$ are in the same component. Otherwise, there exist two nonadjacent vertices of $V_i$ which are not total-monochromatically connected, a contradiction. Thus, the components $C_1,\ldots,C_s$ form a partition of the vertex classes of $G$. Let $C$ be a component of $C_1,\ldots,C_s$. If there is exactly one nontrivial color tree in $C$, it can not be a star. Otherwise, there exist two nonadjacent vertices of the vertex class containing the center, which are not total-monochromatically connected, a contradiction. Hence, there exist at least two internal vertices. Then the nontrivial color tree of $C$ wastes at least $|V(C)|-2+2=|V(C)|$ colors. Suppose $C$ contains $k_c$ ($\geq 2$) nontrivial color trees, denoted by $T_1,\ldots,T_{k_c}$ without loss of generality. If $q_i=1$ for some $i\in\{1,2,\ldots,k_c\}$, then $T_i$ is a star and the center of $T_i$ must be in at least one other nontrivial color tree of $C$ since the vertices of the vertex-class containing the center must be total-monochromatically connected. So we have that
\[\begin{array}{llllll}
\sum_{i=1}^{k_c}(t_i-2+q_i)&\geq\sum_{i=1,q_i\geq 2}^{k_c}(t_i-2+2)+ \sum_{i=1,q_i=1}^{k_c}(t_i-2+1) \\[4pt]
&=\sum_{i=1,q_i\geq 2}^{k_c} t_i+\sum_{i=1,q_i=1}^{k_c}(t_i-1) \\[4pt]
&=\sum_{i=1}^{k_c}t_i-\sum_{i=1,q_i=1}^{k_c}1 \\[4pt]
&\geq |V(C)|+\sum_{i=1,q_i=1}^{k_c}1-\sum_{i=1,q_i=1}^{k_c}1 \\[4pt]
&=|V(C)|.
\end{array}
\]
Then the nontrivial color trees of $C$ waste at least $|V(C)|$ colors. Thus for $1\leq i\leq s$ the nontrivial color trees of $C_i$ waste at least $|V(C_i)|$ colors.
Then $f$ wastes at least $\sum_{i=1}^{s}|V(C_i)|=n$ colors and so $tmc(G)\leq m$. The proof is thus complete.
\qed

\begin{lem}\label{lem5} Let $n$ and $p$ be two integers with $0\leq p\leq n-3$. Then every connected graph $G$ with $n$ vertices and $m=\binom{n}{2}-p$ edges satisfies that $tmc(G)\geq\binom{n}{2}+n-3p$ if $0\leq p\leq \frac{n}{2}$ and $tmc(G)\geq\binom{n}{2}-p$ if $\frac{n}{2}<p\leq n-3$.
\end{lem}

\pf It is trivial for $p=0$, and so assume $1\leq p\leq n-3$. Let $\widetilde{G}$ be the graph obtained from $\overline{G}$ by deleting all the isolated vertices. If $n(\widetilde{G})\leq p+1 \ (\leq n-2)$, then we can find at least two vertices $v_1,v_2$ of degree $n-1$ in $G$. Take a star $S$ with $E(S)=\{v_1v:v\in V(\widetilde{G})\}$. We give all the edges and the internal vertex in $S$ one color, and every other edge and vertex in $G$ a different fresh color. Obviously, it is a TMC-coloring of $G$, which wastes at most $n(\widetilde{G})$ colors. If $n(\widetilde{G})\geq p+2$, say $n(\widetilde{G})=p+t \ (t\geq 2)$, then $\widetilde{G}$ has at least $t$ components since $m(\widetilde{G})=p$. Let $u$ and $v$ be two vertices of $\widetilde{G}$ which are in two different components. We obtain a double star $S'$ in $G$ by connecting $u$ to each vertex in the same component with $v$ of $\widetilde{G}$ and $v$ to the other vertices of $\widetilde{G}$. Assign all the edges and internal vertices in $S'$ one color, and all the other edges and vertices in $G$ different new colors. Clearly, this is a TMC-coloring of $G$, which wastes $n(\widetilde{G})$ colors. If $1\leq p\leq \frac{n}{2}$, then $n(\widetilde{G})\leq 2p$ since $m(\widetilde{G})=p$, implying $tmc(G)\geq m+n-2p=\binom{n}{2}+n-3p$. If $\frac{n}{2}<p\leq n-3$, then we have that $tmc(G)\geq \binom{n}{2}-p$ since $n(\widetilde{G})\leq n$.

The proof is now complete.
\qed

Now we are ready to prove Theorem \ref{thm3}.

{\bf Proof of Theorem 3:} Clearly, $f(n,3)=n-1$, so the assertion holds for $k=3$.

Suppose that $k=\binom{t}{2}+t+2-s$ where $0\leq s\leq t-1$ and $2\leq t\leq n-2$, namely $4\leq k\leq\binom{n}{2}-n+3$. If a connected graph $G$ with $n$ vertices satisfies $m(G)\geq n+k-t-2$, then $l(G)\geq t$ by Lemma \ref{lem3} since $n+k-t-2=n+\binom{t}{2}-s\geq n+\binom{t-1}{2}$. By Theorem \ref{thm1}, we have that  $tmc(G)\geq m-n+2+l(G)\geq n+k-t-2-n+2+t=k$. Thus $f(n,k)\leq n+k-t-2$. To show $f(n,k)\geq n+k-t-2$, it suffices to find a connected graph $G_k$ on $n$ vertices such that $m(G_k)=n+k-t-3$ and $tmc(G_k)<k$.
Take $G_k$ as the graph $G_{t,s}$ described in Lemma \ref{lem0} such that
$m(G_k)=n+\binom{t}{2}-1-s=n+k-t-3$ and $l(G_k)=t$. By Lemma \ref{lem0}, we have that $tmc(G_k)=m(G_k)-n+2+l(G_k)=n+k-t-3-n+2+t=k-1<k$.

Assume that $\binom{n}{2}-n+4\leq k=\binom{n}{2}-q\leq \binom{n}{2}+n-3\lfloor\frac{n}{2}\rfloor$ except for $n$ is odd and $k=\binom{n}{2}+n-3\lfloor\frac{n}{2}\rfloor$.
For a connected graph $G$ with $n$ vertices satisfies $m(G)=\binom{n}{2}-q'\geq k\ (q'\leq q)$, it follows from Lemma \ref{lem5} that $tmc(G)\geq \binom{n}{2}+n-3q'\geq \binom{n}{2}-q'\geq k$ if $0\leq q'\leq\frac{n}{2}$ and $tmc(G)\geq \binom{n}{2}-q'\geq k$ if $\frac{n}{2}<q'\leq q$, implying $f(n,k)\leq k$. To show $f(n,k)\geq k$, it suffices to find a connected graph $G_k$ on $n$ vertices such that $m(G_k)=k-1=\binom{n}{2}-q-1$ and $tmc(G_k)<k=\binom{n}{2}-q$. Take $G_k$ as $G_n^t$ such that $t=2(p+1)-n$ ($p=q+1$) and $m(G_k)=\binom{n}{2}-q-1$. By Lemma \ref{lem4}, we have that $tmc(G_k)=m(G_k)=k-1<k$.

Suppose that $\binom{n}{2}+n-3(r+1)<k\leq\binom{n}{2}+n-3r\ (0\leq r\leq \lfloor\frac{n}{2}\rfloor-1$) or $n$
is odd, $r=\lfloor\frac{n}{2}\rfloor$ and $k=\binom{n}{2}+n-3\lfloor\frac{n}{2}\rfloor$. If a connected graph $G$ on $n$ vertices satisfies $m(G)=\binom{n}{2}-r'\geq \binom{n}{2}-r\ (r'\leq r)$, then $tmc(G)\geq \binom{n}{2}+n-3r'\geq \binom{n}{2}+n-3r\geq k$ by Lemma \ref{lem5}. Thus, $f(n,k)\leq \binom{n}{2}-r$. To show $f(n,k)\geq \binom{n}{2}-r$, it suffices to find a connected graph $G_k$ on $n$ vertices such that $m(G_k)=\binom{n}{2}-r-1$ and $tmc(G_k)<k$. For the case that $\binom{n}{2}+n-3(r+1)<k\leq\binom{n}{2}+n-3r$, where $0\leq r\leq \lfloor\frac{n}{2}\rfloor-1$, take $G_k$ as a complete multipartite graph $K_{n_1,\ldots,n_{n-(r+1)}}$ with $n_1 =\ldots =n_{r+1}= 2$ and $n_{r+2}=\ldots=n_{n-(r+1)}=1$.
It can be checked that $m(G_k)=\binom{n}{2}-r-1$ and $tmc(G_k)=m(G_k)+n-(r+1)-(r+1)=\binom{n}{2}+n-3(r+1)<k$ by Lemma \ref{lem1}. For the case that $n$ is odd, $r=\lfloor\frac{n}{2}\rfloor$ and $k=\binom{n}{2}+n-3\lfloor\frac{n}{2}\rfloor$, take $G_k$ as $G_n^3$ such that $m(G_k)=\binom{n}{2}-\lfloor\frac{n}{2}\rfloor-1$. By Lemma \ref{lem4}, we have that $tmc(G_k)=m(G_k)=\binom{n}{2}-\lfloor\frac{n}{2}\rfloor-1<k$.

The proof is thus complete.\qed

\section{Proof of Theorem \ref{thm4}}

In order to prove Theorem \ref{thm4}, we need the following lemma. Recall that $\binom{1}{2}=0$.

\begin{lem}\label{lem6} Let $G$ be a connected graph with $n$ vertices and $m$ edges. If $\binom{n-t}{2}+t(n-t)\leq m\leq \binom{n-t}{2}+t(n-t)+(t-2)$ for some $t\in\{2,...,n-1\}$, then $tmc(G)\leq m+n-t$. Moreover, the bound is sharp.
\end{lem}

\pf We are given a simple extremal TMC-coloring $f$ of $G$. Since $2\leq t\leq n-1$, we have $m\leq \binom{n}{2}-1$. Then $G$ is not a complete graph and so there is at least one nontrivial color tree. Suppose that $f$ consists of $k$ nontrivial color trees, denoted by $T_1,...,T_k$ where $t_i=|V(T_i)|$ and $q_i=q(T_i)$ for $1\leq i\leq k$. Since $T_i$ has $t_i-1$ edges and $q_i$ internal vertices, it wastes $t_i-2+q_i$ colors. In order to show $tmc(G)\leq m+n-t$, we just need to show that $f$ wastes at least $t$ colors, i.e. $\sum_{i=1}^{k}(t_i-2+q_i)\geq t$. Next it suffices to show that $\sum_{i=1}^{k}(t_i-2)\geq t-1$ since $\sum_{i=1}^{k}q_i\geq 1$. Note that each $T_i$ can total-monochromatically connect at most $\binom{t_i-1}{2}$ pairs of nonadjacent vertices in $G$. Then we have $$\sum_{i=1}^{k}\binom{t_i-1}{2}\geq\binom{n}{2}-m.$$
Suppose $\sum_{i=1}^{k}(t_i-2)<t-1$, that is, $\sum_{i=1}^{k}(t_i-1)<t-1+k$. Since $T_i$ is nontrivial, we have $t_i-1\geq 2$. Then $1\leq k\leq t-2$. By straight forward convexity, the expression $\sum_{i=1}^{k}\binom{t_i-1}{2}$, subject to $t_i-1\geq 2$, is maximized when $k-1$ of the $t_i'$s equal 3 and one of the $t_i'$s, say $t_k$, is as large as it can be, namely, $t_k-1$ is the largest integer smaller than $(t-1+k)-2(k-1)=t-k+1$. Hence $t_k-1=t-k$. Even in this extremal case, we have that  $$\sum_{i=1}^{k}\binom{t_i-1}{2}\leq(k-1)+\binom{t-k}{2}
\leq\binom{t-1}{2}.$$
In fact, $$\binom{t-1}{2}+m\leq\binom{t-1}{2}+\binom{n-t}{2}+t(n-t)+(t-2)
=\binom{n}{2}-1.$$
Hence, $\sum_{i=1}^{k}\binom{t_i-1}{2}\leq \binom{n}{2}-m-1<\binom{n}{2}-m$, a contradiction.

Next we will show that the bound is sharp. Let $G^{*}$ be the graph defined as follows: first take a complete $(n-t+1)$-partite graph with vertex-classes $V_1,...,V_{n-t+1}$ such that $|V_j|=1$ for $1\leq j\leq n-t$ and $V_{n-t+1}=t$; then add the remaining (at most $t-2$) edges to $V_{n-t+1}$ randomly. Clearly, $G^*$ has a spanning subgraph isomorphic to a complete $(n-t+1)$-partite graph
$K_{1,\ldots,1,t}$. By Proposition \ref{pro1} and Theorem \ref{thm1}, it follows that $tmc(G)\geq m+n-t$. Hence, $tmc(G)=m+n-t$.\\
\qed

\noindent$\displaystyle$\textbf{Proof of Theorem 4.} It is trivial for the case that $k=\binom{n}{2}+n$. If $k=\binom{n}{2}+n-1$, we have $g(n,k)\leq \binom{n}{2}-1 $ since $tmc(G)=\binom{n}{2}+n$ for a complete graph $G$. If a connected graph $G$ on $n$ vertices satisfies $m(G)\leq \binom{n}{2}-1$, then there exist two nonadjacent vertices which are total-monochromatically connected by a nontrivial color tree and so it wastes at least two colors. implying $tmc(G)\leq \binom{n}{2}+n-3<k$. Thus, $g(n,k)\geq\binom{n}{2}-1$ and so $g(n,k)=\binom{n}{2}-1$.

For $\binom{n-t}{2}+t(n-t-1)+n\leq k\leq \binom{n-t}{2}+t(n-t)+n-2$ where $2\leq t\leq n-1$, if a connected graph $G$ on $n$ vertices satisfies $m(G)\leq k-n+t(\leq\binom{n-t}{2}+t(n-t)+t-2)$, then $tmc(G)\leq m(G)+n-t\leq k$ by Lemma \ref{lem6}. Hence, $g(n,k)\geq k-n+t$. To show  $g(n,k)\leq k-n+t$, it suffices to find a connected graph $G$ on $n$ vertices such that $m(G)=k-n+t+1$ and $tmc(G)>k$. If $t=2$, then $k=\binom{n}{2}+n-3$ and take $G$ as a complete graph $K_n$. Hence $tmc(G)=\binom{n}{2}+n=k+3>k$. If $t\geq 3$, then take $G$ as the graph $G^{*}$ described in Lemma \ref{lem6} such that $m(G)=k-n+t+1$. It follows from Lemma \ref{lem6} that $tmc(G)=m(G)+n-t=k+1>k$ for $\binom{n-t}{2}+t(n-t-1)+n\leq k\leq \binom{n-t}{2}+t(n-t)+n-3$, and $tmc(G)=m(G)+n-(t-1)=k+2>k$ for $k=\binom{n-t}{2}+t(n-t)+n-2$. Thus $g(n,k)=k-n+t$.

For $k=\binom{n-t}{2}+t(n-t)+n-1$ where $2\leq t\leq n-1$, if a connected graph $G$ on $n$ vertices satisfies $m(G)\leq k-n+t-1(=\binom{n-t}{2}+t(n-t)+t-2)$, then $tmc(G)\leq m(G)+n-t\leq k-1<k$ by Lemma \ref{lem6}. Hence, $g(n,k)\geq k-n+t-1$. To show  $g(n,k)\leq k-n+t-1$, it suffices to find a connected graph $G$ on $n$ vertices such that $m(G)=k-n+t$ and $tmc(G)>k$. If $t=2$, take $G$ as the complete graph $K_n$ and then $tmc(G)=\binom{n}{2}+n=k+2>k$. If $t\geq 3$, take $G$ as the graph $G^{*}$ described in Lemma \ref{lem6} such that $m(G)=k-n+t$. It follows from Lemma \ref{lem6} that $tmc(G)=m(G)+n-(t-1)=k+1>k$. Thus, $g(n,k)=k-n+t-1$.

The proof is now complete.
\qed

\end{document}